\documentclass{amsart}

\newtheorem{theorem}{Theorem}[section]
\newtheorem{lem}[theorem]{Lemma}
\newtheorem{pro}[theorem]{Proposition}

\theoremstyle{definition}

\theoremstyle{remark}
\newtheorem{rem}[theorem]{Remark}

\numberwithin{equation}{section}

\def\G{\Gamma}
\def\R{\mathbb{R}}
\def\O{\Omega}
\def\proof{{\bf Proof\,:\,\,}}
\def\a{\alpha}
\def\e{\epsilon}
\def\g{\gamma}

\begin{document}
\title []{Weighted a priori estimates for the solution of the homogeneous Dirichlet problem for the powers of the Laplacian Operator .}

\author{Ricardo G. Dur\'an}
\address{Departamento de Matem\'atica, Facultad de Ciencias Exactas y Naturales,
Universidad de Buenos Aires, 1428 Buenos Aires, Argentina} \email{rduran@dm.uba.ar}
\author{Marcela Sanmartino}
\address{Departamento de Matem\'atica, Facultad de Ciencias Exactas,
Universidad Nacional de La Plata, 1900 La Plata (Buenos Aires), Argentina }

\email{tatu@mate.unlp.edu.ar}
\author{ Marisa Toschi}
\address{Departamento de Matem\'atica, Facultad de Ciencias Exactas,
Universidad Nacional de La Plata, 1900 La Plata (Buenos Aires), Argentina }

\email{mtoschi@mate.unlp.edu.ar}



\begin{abstract}
Let $u$ be a weak solution of $ (-\Delta)^m u=f $ with Dirichlet
boundary conditions in a smooth bounded domain $\Omega \subset
\mathbb{R}^n$.

 Then, the main goal of this paper is to prove the
following a priori estimate:
$$
\|u\|_{W^{2m,p}_\omega(\Omega)} \le C\,
\|f\|_{L^p_\omega(\Omega)},
$$
where $\omega$ is a weight in the Muckenhoupt class $A_p.$
\end{abstract}

\maketitle
 \baselineskip 8 mm

\section{Introduction}

We will use the standard notation for Sobolev spaces and for
derivatives, namely, if $\alpha$ is a multi-index,
$\alpha=(\alpha_1, \alpha_2, \dots,\alpha_n)\in {\rm Z}\!\!{\rm
Z}^n_+$ we denote $|\alpha|=\sum_{j=1}^n \alpha_j,\
 D^{\alpha}=\partial_{x_1}^{\alpha_1}...\partial_{x_n}^{\alpha_n}$
and
$$
W^{k,p}(\O)=\{v\in L^p(\O)\,:\, D^\alpha v\in L^p(\O) \quad
\forall \, |\alpha|\le k\}.
$$
For $u\in W^{k,p}(\O)$, its norm is given by
$$\|u\|_{W^{k,p}(\O)}=\sum_{|\alpha|\leq k}\|D^\alpha u\|_{L^p(\O)}\, .$$
We consider the homogeneous problem
\begin{eqnarray}\label{dp}\left\{\begin{array}{ccc}
(-\Delta)^m u=f &\mbox{ in }\Omega\\  \\
\left(\frac{\partial}{\partial \nu}\right)^{j}u=0 &\mbox{ in
}\partial\Omega & 0\leq j\leq m-1,
\end{array}\right.
\end{eqnarray}
where $\frac{\partial}{\partial \nu}$ is the normal derivative.

In the classic paper \cite{ADN}, the authors obtained a priori
estimates for solutions of (\ref{dp}) for smooth domain $\Omega$
given by
$$
\|u\|_{W^{2m,p}(\Omega)} \le C\, \|f\|_{L^p(\Omega)}\, .
$$
Where a key tool to prove those estimates was the
Calder\'on-Zygmund theory for singular integral operators.

On the other hand, after the pioneering work of Muckenhoupt
\cite{Mu}, a lot of work on continuity in weighted norms has been
developed. In particular, weighted estimates for a wide class of
singular integral operators has been obtained for weights in the
class of Muckenhoupt $A_p$. Therefore, it is a natural question
whether analogous weighted a priori estimates can be proved for
the derivatives of solutions of elliptic equations.

For the Laplace equation ($m=1$), it was proved in \cite{DST} that
for a weight $\omega$ belonging to the Muckenhoupt class $A_p$
$$
\|u\|_{W^{2,p}_\omega(\Omega)} \le C\, \|f\|_{L^p_\omega(\Omega)}
$$
on a bounded domain $\Omega$ with $\partial\O\in C^2$.

The goal of this paper is to extend the results of \cite{DST} for
powers of the Laplacian operator with homogeneous Dirichlet
boundary conditions, i.e. it is to prove that
\begin{equation}\label{main}
\|u\|_{W^{2m,p}_\omega(\Omega)} \le C\,
\|f\|_{L^p_\omega(\Omega)},
\end{equation}
for $\omega\in A_p$, where the constant $C$ depends on $\O$, $m$,
$n$ and the weight $\omega$.

The main ideas for the proof of these estimates are similar to
those given in \cite{DST}. However, non trivial technical
modifications are needed because, for  $m\geq 2$, the Green
function is not positive in general and therefore, we cannot apply
the maximum principle.

\section{Preliminaries}
\setcounter{equation}{0}

 We consider the problem (\ref{dp}) in a bounded domain $\O$ with
 $\partial\O\in C^{6m+4}$ for $n=2$ and
$\partial\O\in C^{5m+2}$ for $n> 2$
(the regularity on the
boundary is necessary in order to use the results of the Green
function given in  \cite{K}).

The solution of (\ref{dp}) is given by
\begin{equation}
\label{u1} u(x)=\int_\O G_m(x,y)\, f(y)\, dy
\end{equation}
where $G_m(x,y)$ is the Green function of the operator
$(-\Delta)^m$ in $\Omega$ which can be written as
\begin{equation}
\label{u2} G_m(x,y)=\G(x-y)+ h(x,y)
\end{equation}
where $\G(x-y)$ is a fundamental solution and $h(x,y)$ satisfies
\begin{equation*}
\label{h}\left\{\begin{array}{ccc}(-\Delta_x)^m h(x,y)=0  &
x\in \Omega \\
\\ \left(\frac{\partial}{\partial \nu}\right)^{j}h(x,y)
=-\left(\frac{\partial}{\partial\nu}\right)^{j}\G(x-y) &
x\in\partial\Omega & 0\leq j\leq m-1
 \end{array} \right.
 \end{equation*}
  for each fixed $y\in \Omega$. \  \\
Then
\begin{equation}\label{exph}
h(x,y)=-\sum_{j=0}^{m-1}\int_{\partial \O}K_j(y,P)\,
\left(\frac{\partial}{\partial\nu}\right)^j\G(P-x)\, dS
\end{equation}
where $K_j(y,P)$ are the Poisson kernels and $dS$ denotes the
surface measure on $\partial \O$.
\ \\

We recall that any fundamental solution associated to (\ref{dp}) is smooth away from the
origin and it is homogeneous of degree $2m-n$ if $n$ is odd or if $ 2m<n$
and the logarithmic function appears if $n$ is even and $2m\geq n$.
However, in both cases we have the known estimates of the Green
function $G_m(x,y)$ and the Poisson kernels $K_j(x,y)$. In what follows the letter $C$ will denote a generic constant not necessarily the same at each occurrence.
\begin{equation}\label{k_1}
|D_x^\alpha G_m(x,y)|\leq C\ \ \ \ \mbox{for}\ |\alpha|<2m-n,
\end{equation}
\begin{equation}\label{k_2}
|D_x^\alpha G_m(x,y)|\leq C\, \log\left(\frac{2\, diam(\Omega)}{|x-y|}\right)\ \ \ \ \mbox{for}\ |\alpha|=2m-n,
\end{equation}
\begin{equation}\label{k_3}
|D_x^\alpha G_m(x,y)|\leq C\, |x-y|^{2m-n-|\alpha|} \ \ \ \ \mbox{for}\ |\alpha|>2m-n,
\end{equation}
\begin{equation}\label{g_1}
 |D_{x}^\alpha G_m(x,y)|\leq C\, \frac{1}{|x-y|^{n}}\ \min\left\{1,\,
\frac{d(y)}{|x-y|}\right\}^m\ \ \ \ \mbox{for}\ |\alpha|=2m,
\end{equation}
 \begin{equation}\label{g_2}
|K_j(x,y)|\leq C\, \frac{d(x)^{m}}{|x-y|^{n-j+m-1}}\ \ \ \ \mbox{for}\ 0\leq j\leq m-1,
\end{equation}
where $d(x):= \mbox{dist}(x,\partial\Omega)$ (see \cite{K} for (\ref{k_1}), (\ref{k_2}) and (\ref{k_3}) and \cite{D-S} for (\ref{g_1}) and (\ref{g_2})).

\section{The estimates for the derivatives of $u$}
\setcounter{equation}{0}

In this section we state pointwise estimates for the first $2m-1$ derivatives
of the function $u$ and a weak estimate  for the $2m$ derivative.
These estimates will be allow to proof the main result of this work.

\begin{lem}\label{me2m}
Let $u(x)$ be solution of the problem (\ref{dp}). Then, for  $|\alpha|\leq 2m-1$ we have
\begin{eqnarray*} |D^\alpha_x
u(x)|\leq C\, M f(x),
\end{eqnarray*}
where $Mf(x)$ is the usual Hardy- Littlewood maximal function
of $f$.
\end{lem}

\hspace{-0.45cm}\proof
\begin{eqnarray*}\label{ma2m-n}
|D^\alpha_x u(x)|&\leq&\int_{\O}|D^\alpha G_m(x,y)|\, |f(y)|\, dy
\\&&\\
&\leq& C\, \int_{\O} \frac{|f(y)|}{|x-y|^{n-1}}\, dy
\leq C\,  M f(x),
\end{eqnarray*}
by (\ref{k_1}), if $2m-n+1\leq |\alpha|\leq
2m-1$ and by (\ref{k_2}) and (\ref{k_3}), if $|\alpha|\leq 2m-n$.
 \begin{flushright}
$\square$
  \end{flushright}

 \begin{pro}\label{green}
  Given  two measurable functions $f$ and $g$ in $\O$, for $|\alpha|= 2m$ we have that
\begin{eqnarray*} \label{08}\int_{D}|D_x^\alpha G_m(x,y)\, f(y)\, g(x)|\,
 dy\, dx\leq C\,\left( \int_\Omega Mf(x)\, |g(x)|\, dx +  \int_\Omega Mg(y)\, |f(y)|\,
 dy\right),
 \end{eqnarray*}
where $ D:=\{(x,y)\in \Omega\times\Omega\, :\,
|x-y|>d(x)\}$.
  \end{pro}

\hspace{-0.45cm}\proof We write $D=D_1\cup D_2$, where
 \begin{eqnarray*}D_1=\{(x,y)\in D : d(y)\leq 2\, d(x)\}
& \mbox{ and } & D_2=\{(x,y)\in D : d(y)> 2\, d(x)\}.
\end{eqnarray*}
Then, using (\ref{g_1}) we have
\begin{eqnarray}\label{m}\nonumber\int_{D}| D_{x}^\alpha G_m(x,y)\, f(y)\, g(x)|\,
dy\, dx
&\leq&
  \int _{D}\frac{d(y)^m}{|x-y|^{n+m}}\, |f(y)|\, |g(x)|\, dy\, dx
\\\nonumber&&\\
\nonumber&\leq&
 2^m\, \int _{D_1}\frac{d(x)^m}{|x-y|^{n+m}}\, |f(y)|\, |g(x)|\, dy\, dx
\\\nonumber&&\\
&+& \int _{D_2}\frac{d(y)^m}{|x-y|^{n+m}}\, |f(y)|\, |g(x)|\, dy\,
dx=I+II.
\end{eqnarray}

Calling $\O_k(x)= \{z\in\O \, : 2^k d(x)\leq |x-z|<2^{k+1} d(x)\}$,
\begin{eqnarray*}
\int _{D_1}\frac{d(x)^m}{|x-y|^{n+m}}\, |f(y)|\,
|g(x)|\, dy\, dx
 &\leq&
  \int_{\Omega} \sum _{k=1}^\infty  \int _{\O_k(x)}
  \frac{d(x)}{|x-y|^{n+1}}\,
|f(y)|\, dy\, |g(x)|\, dx
\\ &&\\
 &=&
  \int_\Omega A(x)\, |g(x)|\, dx
\end{eqnarray*}
with
\begin{eqnarray*}
 A(x)&\leq&\sum_{k=1}^\infty
  \int_{\{|x-y|< 2^{k+1}d(x)\}}
  \frac{d(x)}{|x-y|^{n+1}}\,
|f(y)|\, dy \leq 2^n \sum _{k=1}^\infty \frac{1}{2^k}\,
  Mf(x)=2^n\, Mf(x).
  \end{eqnarray*}

In order to estimate the term II in (\ref{m}), we first observe
that for $(x,y)\in D_2$, we have that $|x-y|\geq \frac{1}{2}\, d(y)$. Then

\begin{eqnarray*}\label{D_2}
\int _{D_2}\frac{d(y)^m}{|x-y|^{n+m}}\, |f(y)|\,
|g(x)|\, dy\, dx
&\leq&
  C\int_{\Omega} \sum _{k=1}^\infty
  \int _{\O_{k-1}(y)}
  \frac{d(y)}{|x-y|^{n+1}}\,
|g(x)|\, dx\, |f(y)|\, dy
\\&&\\
 &=&\int_\Omega B(y)\, |f(y)|\, dy
\end{eqnarray*}
and therefore, by the same arguments used before we have that
$$B(y)\leq 2^{n+1}\, Mg(y)$$
and the Proposition is proved.
 \begin{flushright}
 $\square$
 \end{flushright}

 In order to see how to estimate in $\Omega\setminus D$, we consider separately
 the function $h$ and $\Gamma$ involved in $G_m.$

 \begin{pro} \label{acoth}
If $|\alpha|\geq 2m-n+1$, there exists a constant $C$ such that

\begin{eqnarray}
\label{02} |D^\a h(x,y)| \le C \,d(x)^{2m-n-|\a|}
\end{eqnarray}
  for $|x-y|\leq d(x).$
\end{pro}

\hspace{-0.45cm}\proof
  In view of (\ref{exph}) we must find estimates for
  $D_x^\alpha(\frac{\partial}{\partial\nu})^j\G(P-x)$ and
 $ K_j(y,P)$.

From the general properties of the fundamental solution
$\Gamma(x-y)$ we have that
\begin{equation}\label{fundamental}
\left|D^\alpha_x(\frac{\partial}{\partial\nu})^j\G(P-x)\right|\leq
C\, |P-x|^{2m-n-|\alpha|-j}
\end{equation}
for $|\alpha|+j\geq 2m-n+1$, and for $0\leq j\leq m-1$, by
  (\ref{g_2}) we have that
  \begin{equation}\label{poisson}
 |K_j(y,P)|\leq C\, \frac{d(y)^m}{|y-P|^{n-j+m-1}}
 \end{equation}
 for $y\in\Omega$ and $P\in \partial\Omega$.

Then by (\ref{fundamental}), (\ref{poisson}) and the fact that if
$|x-y|\leq d(x)$ then $d(y)<2\, d(x)$, we have for $|\alpha|+j\geq
2m-n+1$
 \begin{eqnarray*}
|D^\alpha_x h(x,y)| &\leq& C\, \sum_{j=0}^{m-1}\int_{\partial
\Omega} \frac{d(y)^m}{|y-P|^{n-1+m-j}}\,
 |P-x|^{2m-n-|\alpha|-j}\, dS
 \\&&\\
&\leq& C\, d(x)^{2m-n-|\alpha|}\,
\sum_{j=0}^{m-1}\int_{\partial\Omega}
\frac{d(y)^{m-j}}{|y-P|^{n-1+m-j}}\, dS.
 \end{eqnarray*}

  In order to see that each integral is finite we write
 $\partial\O=F_1\cup F_2$, with
 \begin{eqnarray*}F_1=\{P\in \partial\O: |P_0-P|> 2\, d(y)\}
& and& F_2=\{P\in \partial\O: |P_0-P|\leq 2\, d(y)\},
\end{eqnarray*}
where $P_0\in\partial \Omega$ is that $|y-P_0|=d(y).$ And now, the
convergence of these integrals follow in a standard way.
\begin{flushright} $\square$
 \end{flushright}

It follows from the previous Proposition that for each $x\in \O$
and $|\alpha|\geq 2m-n+1$ we have that $D_x^\alpha h(x,y)$ is
bounded uniformly in a neighborhood of $x$ and so
\begin{equation}\label{parte buena}
D_x^\alpha\int_\O h(x,y)\,
f(y)\, dy=\int_\O D_x^\alpha h(x,y)\, f(y)\, dy.
\end{equation}

On the other hand, although $D_x^\alpha\G$ is a singular kernel
for $|\alpha|=2m$, taking $\beta$ such that $|\beta|=2m-1$, we have that
\begin{equation}
\label{int-sing1}
D_{x_i}\int_\O D_x^\beta\G(x-y) \, f(y) \, dy
= Kf(x)+ c(x)f(x)
\end{equation}
where $c$ is a bounded function and $K$ is a Calder\'on - Zygmund operator given by
\begin{equation}
\label{int-sing2}
Kf(x)= \lim_{\e\to 0}K_{\epsilon}f(x), \ \ \mbox{with} \ \ K_{\epsilon}f(x)=\int_{|x-y|>\e} D_x^\alpha\G(x-y)\, f(y)\, dy .
\end{equation}
Here and in what follows we consider $f$ defined in $\R^n$ extending the original $f$ by zero.

Now we are in conditions to give the following estimate:

\begin{theorem}
\label{acoto} Given $g$ a measurable function and $|\alpha|=2m$.
Then there exists a constant $C$ depending only on $n$, $m$ and
$\O$ such that, for any $x\in\O$,
\begin{eqnarray*} \int_\Omega |D^\alpha_{x} u(x)\, g(x)|\, dx &\leq&
 C\, \left(\int_\Omega \widetilde{K}
f(x)\, |g(x)|\, dx + \int_\Omega Mf(x)\, |g(x)|\, dx\right. \\
&&\\ &+&\left.\int_\Omega Mg(y)\, |f(y)|\, dy +\int_\Omega
|f(x)|\,|g(x)|\, dx\right)
\end{eqnarray*}
where $\widetilde{K}f(x)= \sup_{\e >0}\left|K_{\epsilon}f(x)
\right|$.
\end{theorem}

\hspace{-0.45cm}\proof Using the representation formula for $u$,
by (\ref{parte buena}), (\ref{int-sing1}) and (\ref{int-sing2}) we
have that
\begin{eqnarray}
& \displaystyle{D_x^\alpha u(x) =\lim_{\e\to 0} \int_{\e<|x-y|\le
d(x)}D_x^\alpha\G(x-y)\, f(y)\,dy + c(x)f(x)}\nonumber\\\nonumber&&\\
&\displaystyle{+ \int_{|x-y|\le d(x)}D_x^\alpha h(x,y)\, f(y)\, dy
+\int_{|x-y|>d(x)}D_x^\alpha G(x,y)\, f(y)\, dy}\nonumber\\\nonumber&&\\
 &=:
\displaystyle{I + II + III + IV.} \label{D2}
\end{eqnarray}
By the results given above, for $I$, $II$ and $III$ we have
pointwise estimates, and obtain ( in the same way that in
\cite{DST}) that
$$
|I+II+III|\leq C\left(\widetilde K f (x)+ |f(x)| + M\,
f(x)\right).
$$

However, for $IV$ we have just a weak estimate. Indeed, for the
Proposition \ref{green} we have
\begin{eqnarray*}
\int_\O|IV|\, |g(x)|\, dx
 &\leq&
 C\, \left( \int_\Omega Mf(x)\, |g(x)|\, dx + \int_\Omega Mg(y)\, |f(y)|\,
 dy \right)
\end{eqnarray*}
and the Theorem is proved.
\begin{flushright}
 $\square$
 \end{flushright}

\section {Main result}

We can now state and prove our main result. First we recall the definition of the
$A_p$ class for $1<p<\infty$. A non-negative locally integrable function $\omega$ belongs to
$A_p$ if there exists a constant $C$ such that

\begin{eqnarray*}
\left(\frac{1}{|Q|}\, \int_Q \omega(x)\ dx\right)\left(\frac{1}{|Q|}\,
\int_Q
\omega(x)^{-1/(p-1)}\, dx\right)^{p-1}\le C
\end{eqnarray*}
for all cube $Q\subset\mathbb{R}^n$.

For any weight $\omega$, $L^{p}_\omega(\Omega)$ is the space of measurable
functions $f$ defined in $\Omega$ such that
$$
\|f\|_{L_\omega^p(\Omega)}
=\left(\int_\Omega |f(x)|^p\, \omega(x)\, dx
\right)^{1/p}<\infty
$$
and $W^{k,p}_\omega(\Omega)$ is the space of functions such that
$$
\|f\|_{W^{k,p}_\omega(\Omega) }
=\left(\sum_{|\alpha|\leq k} \|D^\alpha f\|_{L^{p}_\omega(\Omega)}^p\right)^{1/p}<\infty.
 $$

\begin{theorem}
\label{mr} Let $\O\subset \R^n$ be a bounded domain such that
$\partial\O$ is of class $C^{6m+4}$ for $n=2$ and $\partial\O$ is of
class $C^{5m+2}$ for $n\geq 2$.
 If $\omega\in A_p$, $f\in L^p_\omega(\Omega)$ and $u$
 a weak solution of (\ref{dp}),
then there exists a constant $C$ depending only on  $n$, $m$,
$\omega$ and $\O$ such that
$$\|u\|_{W^{2m,p}_\omega(\O)}\leq C \,
\|f\|_{L^p_\omega(\O)}.
$$
\end{theorem}

\hspace{-0.45cm}\proof
 Since $M$ is a bounded operator in $L^p_\omega(\O),$ by Lemma \ref{me2m}
 it follows that
 $$\sum_{|\alpha|\leq 2m-1}\|D^\alpha_x u\|_{L^p_\omega(\O)}\leq C\,
 \|f\|_{L^p_\omega(\O)}.$$

 Therefore, it only remains to estimate $\|D^\alpha_x u\|_{L^p_\omega(\O)}$ for $|\alpha|=2m$.

Let $\omega\in A_p$ and $g(x):=(D^\alpha_x u(x))^{p-1}\,
\omega(x)$. By Theorem \ref{acoto} we see that
\begin{eqnarray}\label{acotdenorma}\nonumber
 \int_\O|D_x^\alpha u(x)|^p\,
\omega(x)\, dx &=&\int_\O |D_x^\alpha u(x)|\, g(x)\, dx
\\\nonumber&&\\\nonumber
 &\leq& C\, \left(\int_\O \widetilde{K}f(x)\, |g(x)|\, dx+
 \int_\O Mf(x)\, |g(x)|\, dx\right.
\\\nonumber&&\\
&+& \left.\int_\O Mg(y)\, |f(y)|\, dy+
 \int_\O |f(x)|\, |g(x)|\, dx\right).
 \end{eqnarray}

Since $\tilde{K}$ and $M$ are bounded operators in
$L^p_\omega(\O)$, applying the H\"older inequality, it follows
that
 \begin{eqnarray}\label{marisa}\nonumber
 \int_\Omega \widetilde{K}f(x)\,
|g(x)|\, dx&=& \int_\Omega \widetilde{K}f(x)\, |g(x)|\,
\frac{1}{\omega(x)^{1/p}}\, \omega(x)^{1/p}\, dx\\
\nonumber&&\\\nonumber
 &\leq& \left(\int_\Omega
\widetilde{K}f(x)^p\, \omega(x)\, dx\right)^{1/p}\,
\left(\int_\Omega |g(x)|^q\, \frac{1}{\omega(x)^{q/p}}\, dx
\right)^{1/q}\\\nonumber &&\\
 &\leq&\|f\|_{L^p_\omega(\Omega)}\, \left(\int_\Omega |g(x)|^q\,
\frac{1}{\omega(x)^{q/p}}\, dx\right)^{1/q},
 \end{eqnarray}
 where $\frac{1}{p}+\frac{1}{q}=1$.

In the same way, we obtain that
 \begin{eqnarray}\label{11}
  \int_\Omega Mf(x)\, |g(x)|\, dx
     \leq \|f\|_{L^p_\omega(\Omega)}\, \left(\int_\Omega |g(x)|^q\,
\frac{1}{\omega(x)^{q/p}}\, dx\right)^{1/q}
   \end{eqnarray}
and
   \begin{eqnarray}\label{12}
 \int_\Omega |f(x)|\, |g(x)|\, dx
   &\leq& \|f\|_{L^p_\omega(\Omega)}\, \left(\int_\Omega |g(x)|^q\,
\frac{1}{\omega(x)^{q/p}}\, dx\right)^{1/q}.
   \end{eqnarray}

For the last term in (\ref{acotdenorma}), taking into account that
 $\omega^{-q/p}\in A_q$, we have that
\begin{eqnarray}\label{aq}
 \int_\Omega M g(y)\, |f(y)|\, dy
&\leq & \|f\|_{L^p_\omega(\Omega)}\, \left(\int_\Omega M g(y)^q\,
\frac{1}{\omega(y)^{q/p}}\, dy\right)^{1/q}
 \\\nonumber&&\\\nonumber
&\leq& \|f\|_{L^p_\omega(\Omega)}\, \left(\int_\Omega |g(x)|^q\,
\frac{1}{\omega(x)^{q/p}}\, dx\right)^{1/q}.
 \end{eqnarray}

Then, by (\ref{marisa}), (\ref{11}), (\ref{12}) and (\ref{aq})we
have
\begin{eqnarray}\label{19}\nonumber
\|D_x^\alpha u\|^p_{L^p_\omega(\O)}&\leq&
 C\, \|f\|_{L^p_\omega(\O)}\, \left(\int_\O |g(x)|^q\,
\frac{1}{\omega(x)^{q/p}}\, dx \right)^{1/q}.
\end{eqnarray}

By the definition of $g(x)$,
\begin{eqnarray*}\label{20}\nonumber\left(\int_\Omega |g(x)|^q\,
\frac{1}{\omega(x)^{q/p}}\, dx \right)^{1/q}
 &=&\left(\int_\Omega |D_x^\alpha u|^{(p-1)q}
 \, \omega(x)^q\, \frac{1}{\omega(x)^{q/p}}\, dx\right)^{1/q}\\\nonumber
 &&\\ &=&\left(\int_\Omega |D_x^\alpha u|^{p}\, \omega(x)\, dx\right)^{1/q}
=\|D_x^\alpha u\|_{L^p_\omega(\Omega)}^{p/q}.
 \end{eqnarray*}

Then we obtain
  \begin{equation}\label{siuinW}
 \|D^\alpha u\|_{L^p_\omega(\Omega)}^p \leq C\, \|f\|_{L^p_\omega(\Omega)}\, \|D^\alpha
u\|_{L^p_\omega(\Omega)}^{p/q}
 \end{equation}
 and the Theorem is proved for $u\in W^{2m,p}_\omega(\Omega).$
\  \\

 Finally, we will show that the weak solutio $u$ of (\ref{dp}) belong to $W^{2m,p}_\omega(\Omega):$

We have that $(-\Delta)^m  u=f$, with $f\in L^p_\omega(\O)$, then
there exists a sequence $f_k\in C^\infty (\R^n)$ such that
$\displaystyle{\lim_{k\rightarrow \infty}f_k=f}$ in
$L^p_\omega(\Omega)$ \cite{Chua}.

For each $k$, there exists $u_k\in C^\infty(\O)$ satisfying
\begin{equation*}\left\{\begin{array}{ccc}(-\Delta)^m  u_k=f_k &\mbox{ in }\O
\\  \\ \left(\frac{\partial}{\partial \nu}\right)^{j}u_k=0
&\mbox{ in }\partial\O& 0\leq j\leq m-1. \end{array}\right.
 \end{equation*}

It is easily to see, from Lemma \ref{me2m} that $u_k\in
W_{\omega}^{2m-1,p}(\O),$ and obviously $u_k\in W_{\omega,\,
loc}^{2m,p}(\O)$. Moreover for all compact set $K\subset \O $, we
have
$$\|u_k\|_{W^{2m,p}_\omega(K)}\leq C(K),$$
where $C(K)$ is a constant depending on the measure of $K$.
Indeed, taking $v_k=u_k \varphi$ with $\varphi \in
 C_0^\infty(K)$, it follows that $v_k \in {W^{2m,p}_\omega(\O)}$,
  satisfies (\ref{dp}) with $f= g_k \in  L^p_\omega(\O)$, and we can use
  (\ref{siuinW}).

Then, it follows by the dominated convergence theorem that $u_k\in
W_{\omega}^{2m,p}(\O)$ and  applying (\ref{siuinW}), we have that
$$\|u_k\|_{W^{2m,p}_\omega(\O)}\leq C\, \|f_k\|_{L^p_\omega(\O).}$$

Therefore, $\{u_k\}$ is a Cauchy sequence in $W^{2m,p}_\omega(\O)$
and there exists $v\in W^{2m,p}_\omega(\O)$ such that
$\displaystyle{\lim_{k\rightarrow \infty}u_k=v}$ in $
W^{2m,p}_\omega(\Omega)$. Let see now that $v$ solves (\ref{dp}).

Obviously, $\displaystyle{f=\lim_{k\rightarrow
\infty}f_k=\lim_{k\rightarrow\infty} (-\Delta)^m u_k=(-\Delta)^m
v}$ in $L^p_\omega(\O)$ and by the classical trace theorems in
Sobolev spaces and the definition of $\omega\in A_p$, it follows
that $v$ satisfies the homogeneous boundary conditions and by
uniqueness of the solution, the Theorem is proved.
 \begin{flushright}
$\square$
 \end{flushright}

\begin{rem}
The result of Theorem \ref{mr} is valid also for $u$ a weak
solution of
\begin{equation*}\left\{\begin{array}{ccc}\mathcal{L}  u=f &\mbox{ in }\O
\\  \\ \mathcal{B}_j u=0
&\mbox{ in }\partial\O& 0\leq j\leq m-1 \end{array}\right.
 \end{equation*}
when $\mathcal{L}\, :=\, \sum_{|\alpha|\, \leq\, 2m}\,
a_{\alpha}\, D^{\alpha}$ is uniformly elliptic and
$\mathcal{B}_j\, :=\, \sum_{\, |\alpha|\, \leq\, j}\, b_{\alpha}\,
D^{\alpha},$ $ 0\, \leq j \leq m-1$ are the boundary operators
defined in \cite{ADN}.

Indeed, we define $l_1> \max_{j} (2m - j)$ and $l_0 = \max_j (2m -
j)$. If the coefficients $a_\alpha\in
C^{l_1+1}(\overline{\Omega})$, $b_j\in C^{l_1+1}(\partial\Omega)$
and $\partial\Omega\in C^{l_1+2m+1}$
 we have that the Green
function $G_m$ and the Poisson kernels $K_j$ for $0\leq j \leq m -
1$ exist whenever $l_1 > 2 (l_0 + 1)$ for $n = 2$ and $l_1 > \frac
{3} {2}\,  l_0$ for $n \geq 3$.

Moreover, wherever they are defined, the Green function and the
Poisson kernels of the operator $\mathcal{L}$ with these boundary
conditions satisfy the estimates (\ref{k_1}),\ (\ref{k_2}),\
(\ref{k_3}),\ (\ref{g_1}) and (\ref{g_2}) (see \cite{D-S} and
\cite{K}).
\end{rem}

\begin{rem}
Using the fact that $d(x)^\beta \in A_p$ for $-1<\beta < p-1$ and
some imbedding Theorems for weighted Sobolev spaces (see
\cite{DST}) we have as a consequence of the main result
\begin{theorem}
\label{pq} Let $\Omega\subset \mathbb{R}^n$ be a bounded domain as
above, $f\in L^p_{d^{\gamma}}(\Omega)$, with $\g=k\beta$, where
$k\in \mathbb{N}$ and $0\leq\beta\leq 1$. If $u$ be the solution
of problem (\ref{dp}), $0 \leq \gamma < p-1$ and
$\displaystyle{\frac{1}{p}-\frac{1}{q}\leq \frac{2m}{n+k}}$ \, (
with $q<\infty$ when $2mp=n+k$), then there exists a constant $C$
depending only on $\gamma$ , $p$, $q$, $n$ and $\Omega$ such that

\begin{equation}\label{aplication}
\|u\|_{L^{q}_{d^{\gamma}}(\Omega)}\leq C\,
\|f\|_{L^p_{d^{\gamma}}(\Omega)}.
\end{equation}
\end{theorem}
 Finally, as a particular case of (\ref{aplication}) taking $\g=m$ we
 have that

$$
\|u\|_{L^q_{d^m}(\Omega)}\leq C\, \|f\|_{L^p_{d^m}(\Omega)}
$$

for $p>m+1$ and $\displaystyle{\frac{1}{p}-\frac{1}{q}\leq
\frac{2m}{n+1}}$ ( with $q<\infty$ when $2mp=n+m$).

 This result is proved in \cite{D-S} using different arguments
for the case
$\displaystyle{\frac{1}{p}-\frac{1}{q}<\frac{2m}{n+m}}$.

Our results shows that, at least in the case $p>m+1$, the estimate
remains valid when $\displaystyle{\frac{1}{p}-\frac{1}{q}=
\frac{2m}{n+m}}$.

\end{rem}

\end{document}